\documentclass[letterpaper, 10 pt, conference]{ieeeconf}  

\IEEEoverridecommandlockouts                              
\usepackage{xz_style}
\usepackage{cite}
\renewcommand\theequation{\arabic{section}.\arabic{equation}}

\renewcommand\thesection{\arabic{section}}

\title{\LARGE \bf
Revisiting LQR Control from the Perspective of \\ Receding-Horizon Policy Gradient 
}

\author{Xiangyuan Zhang \qquad Tamer Ba\c{s}ar
\thanks{The authors are with the Department of ECE and CSL, University of Illinois at Urbana-Champaign, Urbana, IL 61801, USA. {\tt\small \{xz7, basar1\}@illinois.edu}. Research of the authors was supported in part by the Air Force Office of Scientific Research (AFOSR) through Grant FA9550-19-1-0353. The latest version is updated on Jan. 2024, where we have removed a required condition on $Q_N$ and improved the presentation.}%
}

\begin{document}

\maketitle
\thispagestyle{empty}
\pagestyle{empty}

\begin{abstract}
We revisit in this paper the discrete-time linear quadratic regulator (LQR) problem from the perspective of receding-horizon policy gradient (RHPG), a newly developed model-free learning framework for control applications. We provide a fine-grained sample complexity analysis for RHPG to learn a control policy that is both stabilizing and $\epsilon$-close to the optimal LQR solution, and our algorithm does not require knowing a stabilizing control policy for initialization. Combined with the recent application of RHPG in learning the Kalman filter, we demonstrate the general applicability of RHPG in linear control and estimation with streamlined analyses. 
\end{abstract}

\section{Introduction}
Model-free policy gradient (PG) methods promise a universal end-to-end framework for controller designs. By utilizing input-output data of a black-box simulator, PG methods directly search the prescribed policy space until convergence, agnostic to system models, objective function, and design criteria/constraints. The general applicability of PG methods leads to countless empirical successes in continuous control, but the theoretical understanding of these PG methods is still in its early stage. Stemmed from the convergence theory of PG methods for general reinforcement learning tasks \cite{agarwal2020optimality, zhang2020global}, a recent thrust of research has specialized the analysis for the convergence and sample complexity of PG methods into several linear state-feedback control benchmarks \cite{fazel2018global, malik2020derivative, li2019distributed, hambly2020policy, zhang2021derivative, perdomo2021stabilizing, hu2022towards}. However, incorporating imperfect-state measurements leads to a deficit of most, if not all, favorable landscape properties crucial for PG methods to converge (globally) in the state-feedback settings \cite{zheng2021analysis, hu2022towards}. Even worse, the control designer now faces several challenges unique to control applications: a) convergence might be toward a suboptimal stationary point without system-theoretic interpretations; b) provable stability and robustness guarantee could be lacking; c) convergence depends heavily on the initialization (e.g., the initial policy should be stabilizing), which would be challenging to hand-craft; and d) algorithm could be computationally inefficient. These bottlenecks blur the applicability of model-free PG methods in real-world control scenarios since the price for each of the above disadvantages could be unaffordable.

On the other hand, classic theories provide both elegant analytic solutions and efficient computational means (e.g., Riccati recursions) to a wide range of control problems \cite{anderson1979optimal, anderson1990optimal, bacsar1995h}. They further reveal the intricate structure in various control settings and offer system-theoretical interpretations and guarantees to their characterized solutions. They suggest that, compared to viewing the dynamical system as a black box and studying the properties of PG methods from a (nonconvex) optimization perspective, it is better to incorporate those properties unique to decision and control into the design of learning algorithms.

In this work, we revisit the classical linear quadratic regulator (LQR) problem from the perspective of the newly-developed receding-horizon PG (RHPG) framework \cite{zhang2023learning}, which integrates Bellman’s principle of optimality into the development of a model-free PG framework. First, RHPG approximates infinite-horizon LQR using a finite-horizon problem formulation and further decomposes the finite-horizon problem into a sequence of one-step sub-problems. Second, RHPG solves each sub-problem recursively using model-free PG methods. To accommodate the inevitable computational errors in solving these sub-problems, we establish the generalized principle of optimality that bounds the accumulated bias by controlling the inaccuracies in solving each sub-problem. We characterize the convergence and sample complexity of RHPG in \S\ref{sec:sample} and emphasize that the RHPG algorithm does not require knowing a stabilizing initial control policy \emph{a priori}. 
 
 \subsection{Literature Review}
We mainly compare with \cite{fazel2018global, malik2020derivative} and \cite{lamperski2020computing, perdomo2021stabilizing, zhao2021learning}, where \cite{fazel2018global, malik2020derivative} are the foundational work in applying policy optimization to LQR and \cite{lamperski2020computing, perdomo2021stabilizing, zhao2021learning} remove the assumption on an initial stabilizing point in LQR by adding a discount factor to the objective function as an extra parameter.

In contrast to \cite{fazel2018global, malik2020derivative} that parametrizes LQR as a single nonconvex (constrained) optimization problem over the policy space, we provide a new parametrization and perspective to learning LQR control. The critical difference between the RHPG algorithm and \cite{fazel2018global, malik2020derivative} is that RHPG incorporates existing theories into the design of model-free learning algorithms, which is more than just exploiting them for the convergence analysis. We view this work as an initial step toward the goal of ``\emph{designing control-specific learning algorithms with performance guarantees}'' rather than ``\emph{analyzing existing learning algorithms for control}''. This line of research is motivated by the observation that viewing the dynamical system as a black box and directly searching in the policy space leads to a deficit of most, if not all, favorable landscape properties beyond LQR \cite{zheng2021analysis, hu2022towards}. This implies that the excellent properties in LQR following the parametrization in \cite{fazel2018global, malik2020derivative}, such as coercivity and gradient domination, are rare, problem-dependent, and hard to generalize. However, when the model information is known, existing theories have provided extremely efficient solutions (e.g., Riccati recursions) to these problems that seemed unsolvable in the PO paradigm. Hence, our rationale is that control settings inherently have more structures (that are problem-agnostic) than a black-box system. One should always exploit these structures in \emph{developing learning-based algorithms with performance guarantees}. In our work and \cite{zhang2023learning, zhang2023global}, we have identified causality and the dynamic programming principle as the fundamental properties in all control settings and exploited them in the model-free learning paradigm. As demonstrated in our work and \cite{zhang2023learning, zhang2023global}, the RHPG framework efficiently solves LQR and the seemed-to-be-unsolvable Kalman filtering problem in the PO fashion, which serves as a fundamental benchmark in output feedback control. As a by-product, the RHPG framework also removes all assumptions inappropriate in model-free learning settings and is more consistent with existing theories in control and estimation. Our work and \cite{zhang2023learning, zhang2023global} together lead to a promising path toward the theoretical foundation of model-free learning in partially observable settings and nonlinear control through the lens of RHPG. 
	
In comparison with \cite{lamperski2020computing, perdomo2021stabilizing, zhao2021learning}, we note that the $\gamma$-discounted LQR problems therein are equivalent to the standard non-discounted LQR with system matrices being $\sqrt{\gamma}A$ and $\sqrt{\gamma}B$. Then, for any $\gamma\in (0, 1)$, the set of stabilizing policies with system matrices being $\sqrt{\gamma}A$ and $\sqrt{\gamma}B$ is strictly larger than that of the un-discounted case. Hence, when $\gamma$ is sufficiently small, one can initialize the PG algorithm with an arbitrary control policy. This removes the requirement of knowing a stabilizing policy in advance, but it comes with the price of solving multiple LQRs instead of only one. Moreover, the criterion for increasing $\gamma$ is more complex than our selection rule for $N$ in RHPG. Besides removing the assumption on the initial stabilizing point, it is more important to determine if the results/insights can be further generalized/extended to more complicated control problems, where the critical difference between our work and \cite{lamperski2020computing, perdomo2021stabilizing, zhao2021learning} appears. In \cite{lamperski2020computing, perdomo2021stabilizing, zhao2021learning}, the landscape of discounted LQR is identical to those in \cite{fazel2018global, malik2020derivative} for un-discounted LQR, with essentially scaled versions of system parameters. Thus, the same difficulties reported in \cite{zheng2021analysis, hu2022towards} will appear when considering the output-feedback setting with an additional discount factor. In contrast, RHPG can be directly extended to solve output-feedback problems \cite{zhang2023global}.

\subsection{Notations}\label{sec:notations}
We use $\|X\|$, $\kappa_X$, and $\rho(X)$ to denote, respectively, the spectral norm, condition number, and the spectral radius of a square matrix $X$. If $X$ is symmetric,  we use $X > 0$ and $X\geq 0$ to denote that $X$ is positive definite (pd) and positive semi-definite (psd), respectively.  For a pd matrix $W$ of appropriate dimensions, we define the $W$-induced norm of $X$ as 
\begin{align*}
	\|X\|^2_W := \sup_{z \neq 0} \frac{z^{\top}X^{\top}WXz}{z^{\top}Wz}.
\end{align*}

\section{Preliminaries}
\subsection{Infinite-Horizon LQR}
Consider the discrete-time linear dynamical system\footnote{For extensions to stochastic LQR with i.i.d. additive noises, as well as the setting with an arbitrary (deterministic) initial state, see \S\ref{sec:discussion}.}
\begin{align}\label{eqn:inf_LQR_dynamics}
	x_{t+1} = Ax_t+ Bu_t, 
\end{align}
where $x_t\in \RR^n$ is the state; $u_t \in \RR^{m}$ is the control input; $A \in \RR^{n\times n}$ and $B \in \RR^{n\times m}$ are system matrices unknown to the control designer; and the initial state $x_0\in\RR^n$ is sampled from a zero-mean distribution $\cD$ that satisfies $\Cov(x_0) = \Sigma_0 > 0$. The goal in the LQR problem is to obtain the optimal controller $u_t=\phi_t(x_t)$ that minimizes the cost
\begin{align}\label{eqn:cost}
	 J_{\infty} := \EE_{x_0\sim\cD}\left[\sum_{t=0}^{\infty} \big( x_t^{\top}Qx_t + u^{\top}_tRu_t \big)\right],
\end{align}
where $Q > 0$ and $R > 0$ are symmetric pd weightings chosen by the control designer. For the LQR problem as posed to admit a solution, we require $(A, B)$ to be stabilizable. Note that here $Q > 0$ implies the observability of $(A, Q^{1/2})$. Then, the unique optimal LQR controller is linear state-feedback, i.e., $u^*_t = -K^*x_t$, and $K^* \in \RR^{m\times n}$, which with a slight abuse of terminology we will call optimal control policy, can be computed by
\begin{align}\label{eqn:infinite_lqr_gain}
	K^* = (R + B^{\top}P^*B)^{-1}B^{\top}P^*A,
\end{align}
where $P^*$ is the unique positive definite (pd) solution to the algebraic Riccati equation (ARE)
\begin{align}\label{eqn:ARE}
	P = Q + A^{\top}PA - A^{\top}PB(R+B^{\top}PB)^{-1}B^{\top}PA.
\end{align}
Moreover, the optimal control policy $K^*$ is guaranteed to be stabilizing, i.e., $\rho(A-BK^*)<1$. Therefore, we can parametrize LQR as an optimization problem over the policy space $\RR^{m\times n}$, subject to the stability condition \cite{fazel2018global}:
\begin{align}\label{eqn:LQR}
	&\hspace{-0.2em}\min_{K}~ J_{\infty}(K) = \EE_{x_0\sim\cD}\left[\sum_{t=0}^{\infty}\big(x_t^{\top}(Q+K^{\top}RK)x_t\big)\right]\\
	&\hspace{2em}\text{s.t.} \quad K \in \cK := \{K \mid \rho(A-BK)<1\}. \label{eqn:ck}
\end{align}
Theoretical properties of model-free (zeroth-order) PG methods in solving \eqref{eqn:LQR} have been well understood \cite{fazel2018global, bu2019LQR, malik2020derivative}. In particular, the objective function \eqref{eqn:LQR}, even though nonconvex, is coercive and (globally) gradient dominated \cite{bu2019LQR}. Hence, if an initial control policy $K_0 \in \cK$ is known {\it a priori}, then any descent direction of the objective value (e.g., vanilla PG) suffices to ensure that all the iterates will remain in the interior of $\cK$ while quickly converging toward the unique stationary point. Removing the assumption on $K_0$ (that an initial stabilizing policy can readily be found) has remained an active research topic \cite{lamperski2020computing, perdomo2021stabilizing, zhao2021learning}.

\subsection{Finite-Horizon LQR}
The finite-$N$-horizon version of the LQR problem is also described by the system dynamics \eqref{eqn:inf_LQR_dynamics}, but with the objective function summing up only up to time $t=N$. Similar to \eqref{eqn:LQR}, we can parametrize the finite-horizon LQR problem as $\min_{\{K_t\}}  J\big(\{K_t\}\big)$, where

\vspace{-0.5em}
\small
\begin{align}\label{eqn:finite_LQR}
	\hspace{-0.3em} J\big(\hspace{-0.1em}\{K_t\}\hspace{-0.1em}\big) \hspace{-0.2em}:=\hspace{-0.15em}  \EE_{x_0\sim \cD}\hspace{-0.3em}\left[\hspace{-0.05em}\sum_{t=0}^{N-1}\hspace{-0.1em}x_t^{\top}\hspace{-0.2em}(Q\hspace{-0.1em}+\hspace{-0.1em}K^{\top}_tRK_t)x_t \hspace{-0.1em}+\hspace{-0.1em} x_{\hspace{-0.1em}N}^{\top}Q_{\hspace{-0.1em}N}x_{\hspace{-0.1em}N}\hspace{-0.2em}\right]\hspace{-0.3em},
\end{align}
\normalsize
and $Q_N$ is a symmetric psd terminal-state weighting to be chosen. The unique optimal control policy in the finite-horizon LQR is time-varying and can be computed by 
\begin{align}\label{eqn:finite_lqr_gain}
	K^*_t = (R+B^{\top}P^*_{t+1}B)^{-1}B^{\top}P^*_{t+1}A,
\end{align}
where $P^*_t$, for all $t \in \{0, \cdots, N-1\}$, are generated by the Riccati difference equation (RDE) starting with $P^*_N = Q_N$:
\begin{align}
	P^*_t &= Q + A^{\top}P^*_{t+1}A  \nonumber\\
	&\hspace{1.5em}- A^{\top}P^*_{t+1}B(R+B^{\top}P^*_{t+1}B)^{-1}B^{\top}P^*_{t+1}A. \label{eqn:RDE}
\end{align}
Theoretical properties of zeroth-order PG methods in addressing \eqref{eqn:finite_LQR} have been studied in \cite{zhang2021derivative, hambly2020policy}. Compared to the infinite-horizon setting \eqref{eqn:LQR}, the finite-horizon LQR problem \eqref{eqn:finite_LQR} is also a nonconvex and gradient-dominated problem, but it does not naturally require the stability condition \eqref{eqn:ck}.

\section{Receding-Horizon Policy Gradient}
\subsection{LQR with Dynamic Programming}
It is well known that the solution of the RDE \eqref{eqn:RDE} converges monotonically to the stabilizing solution of the ARE \eqref{eqn:ARE} exponentially \cite{hassibi1999indefinite}. It then readily follows that the sequence of time-varying LQR policies \eqref{eqn:finite_lqr_gain}, denoted as $\{K_t\}_{t\in \{N-1, \cdots, 0\}}$, converges monotonically to the time-invariant LQR policy $K^*$ as $N\to\infty$. Now, we formally present this non-asymptotic convergence result.

\begin{theorem}\label{lemma:finite_approximation}
Let $A_K^*:=A-BK^*$, use $\|\cdot\|_*$ to denote the $P^*$-induced norm, and define
	\begin{align}\label{eqn:N0}
		N_0 = \frac{1}{2}\cdot \frac{\log\big(\frac{\|Q_N-P^*\|_*\cdot\kappa_{P^*}\cdot \|A_K^*\|\cdot\|B\|} {\epsilon\cdot\lambda_{\min}(R)}\big)}{\log\big(\frac{1}{\|A_K^*\|_*}\big)} + 1.
	\end{align}
	where it holds that $\|A_K^*\|_* <1$. Then, for all $N\geq N_0$, the control policy $K^*_{0}$ computed by \eqref{eqn:finite_lqr_gain} satisfies $\|K^*_{0} - K^*\| \leq \epsilon$ for any $\epsilon > 0$.
\end{theorem}\par  
We provide the proof of Theorem \ref{lemma:finite_approximation} in \S\ref{proof:finite}. Theorem \ref{lemma:finite_approximation} demonstrates that if selecting $N\sim \cO(\log(\epsilon^{-1}))$, then solving the finite-horizon LQR will result in a policy $K^*_0$ that is $\epsilon$-close to $K^*$, for any $\epsilon > 0$. Furthermore, if one chooses a small enough $\epsilon$ such that an $\epsilon$-ball centered at $K^*$ lies entirely in $\cK$, then this condition on  $\epsilon$ constitutes a sufficient condition for $K^*_0$ to be stabilizing, i.e., $K^*_0 \in \cK$.


\subsection{Algorithm Design}

\begin{figure*}
\end{figure*}

\begin{algorithm}[t]
  \label{alg:DP}
 \caption{Receding-Horizon Policy Gradient}
  \SetAlgoLined
  \KwIn{horizon $N$, max iterations $\{T_h\}$, smoothing radius $\{r_h\}$, stepsizes $\{\eta_{h}\}$}
  \For{$h = N-1, \cdots, 0$}{ 
Initialize $K_{h, 0}$ arbitrarily (e.g., the convergent policy from the prev. iter. $K_{h+1, T_{h+1}}$ or $0$)\;
\For{$i = 0, \cdots, T_h-1$ }{// sample PG update via a zeroth-order oracle \\
	Sample $K_{h, i}^{+} = K_{h, i} + r_hU$ and $K_{h, i}^{-} = K_{h, i} - r_hU$, where $U$ is uniformly drawn from the surface of a unit sphere, i.e., $\|U\|_F=1$\;
	Sample $x_h \sim \cD$ and simulate two trajectories with policies $K_{h, i}^{+}$ and $K_{h, i}^{-}$, respectively. Compute values $J_h(K_{h, i}^{+})$ and $J_h(K_{h, i}^{-})$\;
	Compute the estimated PG $\tilde{\nabla} J_h(K_{h, i})\hspace{-0.15em}=\hspace{-0.15em} \frac{mn}{2r_h}\big[J_h(K_{h, i}^{+})\hspace{-0.15em}-\hspace{-0.15em}J_h(K_{h, i}^{-})\big]U$ \\
    and update $K_{h, i+1} = K_{h, i} - \eta_{h}\cdot \tilde{\nabla}  J_h(K_{h, i})$\;
  }}
  Return $K_{0, T_{0}}$\;
\end{algorithm}

We propose the RHPG algorithm (cf., Algorithm \ref{alg:DP}), which first selects $N$ by Theorem \ref{lemma:finite_approximation}, and then sequentially decomposes the finite-$N$-horizon LQR problem backward in time. In particular, for every iteration indexed by $h \in \{N-1, \cdots, 0\}$, the RHPG algorithm solves an LQR problem from $t=h$ to $t=N$, where we only optimize for the current policy $K_h$ and fix all the policies $\{K_t\}$ for $t\in\{h+1, \cdots, N-1\}$ to be the convergent solutions generated from earlier iterations. Concretely, for every $h$, the RHPG algorithm solves the following \emph{quadratic} program in $K_h$:
\begin{align}
	\min_{K_h}\  J_h(K_h) &:= \EE_{x_h\sim\cD}\bigg[\sum^{N-1}_{t=h+1} x_t^{\top}\big(Q + (K_t^*)^{\top}RK^*_t \big)x_t \nonumber\\
	&\hspace{-0.4em}+ x_h^{\top}\big(Q +K_h^{\top}RK_h\big)x_h + x_N^{\top}Q_Nx_N\bigg]. \label{eqn:induction_lqr}
\end{align}
Due to the quadratic optimization landscape of \eqref{eqn:induction_lqr} in $K_h$ for every $h$, applying any PG method with an arbitrary finite initial point (e.g., zero) would lead to convergence to the globally optimal solution of  \eqref{eqn:induction_lqr}.

\subsection{Bias of Model-Free Receding-Horizon Control}\label{sec:bias}

\begin{figure*}
	\centering
	\includegraphics[width=0.8\textwidth]{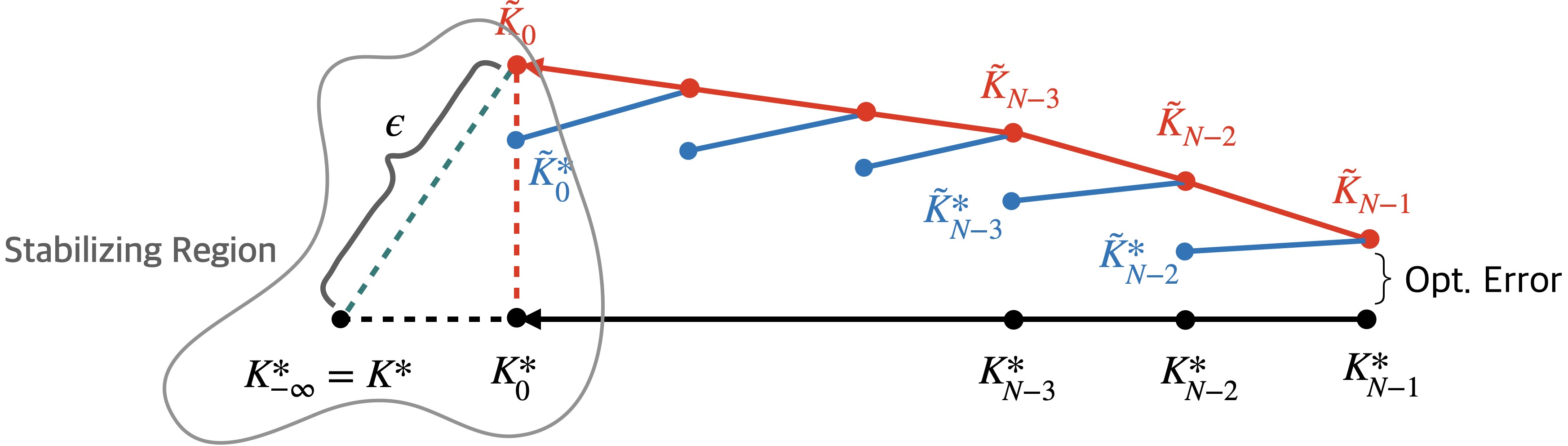}
	\caption{We first show that the output policy $\tilde{K}_{0}$ can be made $\epsilon$-close to $K^*$ in two steps. First, Theorem \ref{lemma:finite_approximation} proves that $K^*_{0}$ is $\epsilon$-close to $K^*$ by selecting $N$ accordingly. Then, Theorem \ref{theorem:LQR_DP} analyzes the backward propagation of the computational errors from solving each subproblem, denoted as $\delta_t:=\tilde{K}_t - \tilde{K}^*_t$ for all $t$, where $\tilde{K}^*_t$ represents the current optimal LQR policy after absorbing errors from all previous iterations. Then, we show that if one requires a small enough optimality gap $\epsilon$ between $\tilde{K}_0$ and $K^*$, then the RHPG output $\tilde{K}_0$ can automatically acquire a closed-loop stability certificate. }\label{fig:proof_sketch}
\end{figure*}

The RHPG algorithm builds on Bellman's principle of optimality, which requires solving each iteration to the exact optimal solution. However, PG methods can only return an $\epsilon$-accurate solution after a finite number of steps. To generalize Bellman's principle of optimality, we analyze how computational errors accumulate and propagate in the (backward) dynamic programming process. In the theorem below, we show that if one solves every iteration of the RHPG algorithm to the $\cO(\epsilon)$-neighborhood of the unique optimum, then the RHPG algorithm will output a policy that is $\epsilon$-close to the infinite-horizon LQR policy $K^*$. 

\begin{theorem}\label{theorem:LQR_DP}
	Choose $N$ according to Theorem \ref{lemma:finite_approximation} and assume that one can compute, for all $h\in\{N-1, \cdots, 0\}$ and some $\epsilon > 0$, a policy $\tilde{K}_h$ that satisfies
	\begin{align*}
		\big\|\tilde{K}_{h}  - \tilde{K}_{h}^*\big\| \sim\cO(\epsilon)\cO(\texttt{poly}(\text{system parameters})),
	\end{align*}
	where $\tilde{K}_{h}^*$ is the optimum of the LQR from $h$ to $N$, after absorbing errors in all previous iterations of Algorithm \ref{alg:DP}. Then, the RHPG algorithm outputs a control policy $\tilde{K}_0$ that satisfies $\big\|\tilde{K}_0 - K^*\big\| \leq \epsilon$. Furthermore, if $\epsilon$ is sufficiently small such that $\epsilon < \frac{1-\|A-BK^*\|_*}{\|B\|}$, then $\tilde{K}_0$ is stabilizing.
\end{theorem}

We illustrate Theorem \ref{theorem:LQR_DP} in Figure \ref{fig:proof_sketch} and defer its proof to \S\ref{proof:LQR_DP}. Additionally, we discuss the implications of Theorem \ref{theorem:LQR_DP} in the following remark.

\begin{remark}
Theorem \ref{theorem:LQR_DP} ensures that if each RHPG iteration satisfies a certain error tolerance level, then the RHPG output $\tilde{K}_0$ will reach an $\epsilon$-neighborhood of $K^*$. Notably, upon selecting a sufficiently small $\epsilon$, such that the $\epsilon$-ball around $K^*$ lies entirely within $\cK$, we can guarantee closed-loop stability of $\tilde{K}_0$ solely based on its near-optimality. This approach differs significantly from existing PG for LQR literature (e.g., \cite{fazel2018global, malik2020derivative}), which requires a stabilizing initial policy and focus on preserving the closed-loop stability during learning. Conversely, RHPG starts with an arbitrary initial point, potentially very far from the stabilizing region (see our numerical experiments in \S\ref{sec:sim}), yet still converges globally to $K^*$. The closed-loop stability certificate then comes for free, given near-optimality.
\end{remark}

Now, it remains to establish the sample complexity for the convergence of (zeroth-order) PG methods in every iteration of the algorithm, which is done next.

\subsection{PG Update and Sample Complexity}\label{sec:sample}
We analyze here the sample complexity of the zeroth-order PG update in solving each iteration of the RHPG algorithm. Specifically, the zeroth-order PG update is defined as
\begin{align}\label{eqn:PG}
	K_{h, i+1} = K_{h, i} - \eta_{h}\cdot \tilde{\nabla}  J_h(K_{h, i})
\end{align}
where $\eta_{h} > 0$ is the stepsize to be determined and $\tilde{\nabla}  J_h(K_{h, i})$ is the estimated PG sampled from a (two-point) zeroth-order oracle. We formally present the sample complexity result in the following proposition.


\begin{proposition}\label{prop:sample}
For all $h\in\{0, \cdots, N-1\}$, choose a constant smoothing radius $r_{h} \sim \cO(\epsilon)$ and a constant stepsize $\eta_{h} \sim \cO(\epsilon^2)$. Then, the zeroth-order PG update \eqref{eqn:PG} converges after $T_h \sim \cO(\frac{1}{\epsilon^2}\log(\frac{1}{\delta\epsilon^2}))$ iterations in the sense that $\big\|K_{h, T_h} - \tilde{K}_{h}^*\big\| \leq \epsilon$ with a probability of at least $1-\delta$.
\end{proposition}

For completeness, we provide a supplementary proof of Proposition \ref{prop:sample} in \S\ref{proof:sample}, which mostly follows existing results in the literature \cite{malik2020derivative}. Combining Theorem \ref{theorem:LQR_DP} with Proposition \ref{prop:sample}, we conclude that if we spend $\tilde{\cO}(\epsilon^{-2}\log(\delta^{-1}))$ iterations in solving every subproblem to $\cO(\epsilon)$-accuracy with a probability of $1-\delta$, for all $h \in \{0, \cdots, N-1\}$, then the RHPG algorithm will output a $\tilde{K}_{0}$ that satisfies  $\|\tilde{K}_0 - K^* \| \leq \epsilon$ with a probability of at least $1-N\delta$. By \eqref{eqn:N0}, this implies that the total iteration complexity of RHPG is also $\tilde{\cO}(\epsilon^{-2}\log(\delta^{-1}))$ with the dependence on various system parameters being polynomial.

We discuss the tradeoffs in selecting $N$ to balance minimizing finite-to-infinite error and minimizing errors through the backward propagation in \S\ref{sec:tradeoff}. To compare our sample complexity bound with the sharpest result in the literature \cite{malik2020derivative}, our dependence on $\epsilon$ matches that of \cite{malik2020derivative}\footnote{Note that the $\tilde{\cO}(\epsilon^{-1})$ sample complexity presented in \cite{malik2020derivative} is for the convergence in objective value (e.g., $f(K) - f(K^*) \leq \epsilon$), and is equivalent to an $\tilde{\cO}(\epsilon^{-2})\cdot\cO(\texttt{poly}(\text{system parameters}))$ sample complexity for the convergence in policy (i.e., $\|K-K^*\| \leq \epsilon$).}. Our sample complexity and that of \cite{malik2020derivative} have polynomial dependence on system parameters. However, it is not clear how to compare the polynomial dependencies between our bounds and those of \cite{malik2020derivative}, and these polynomial factors might affect the overall computational efficiency of both algorithms in a substantial way. We leave this comparison as an important future research topic.

\section{Numerical Experiments}\label{sec:sim}
\begin{figure}
\vspace{0.5em}
	\centering\hspace{-1.5em}
	\includegraphics[width=0.46\textwidth]{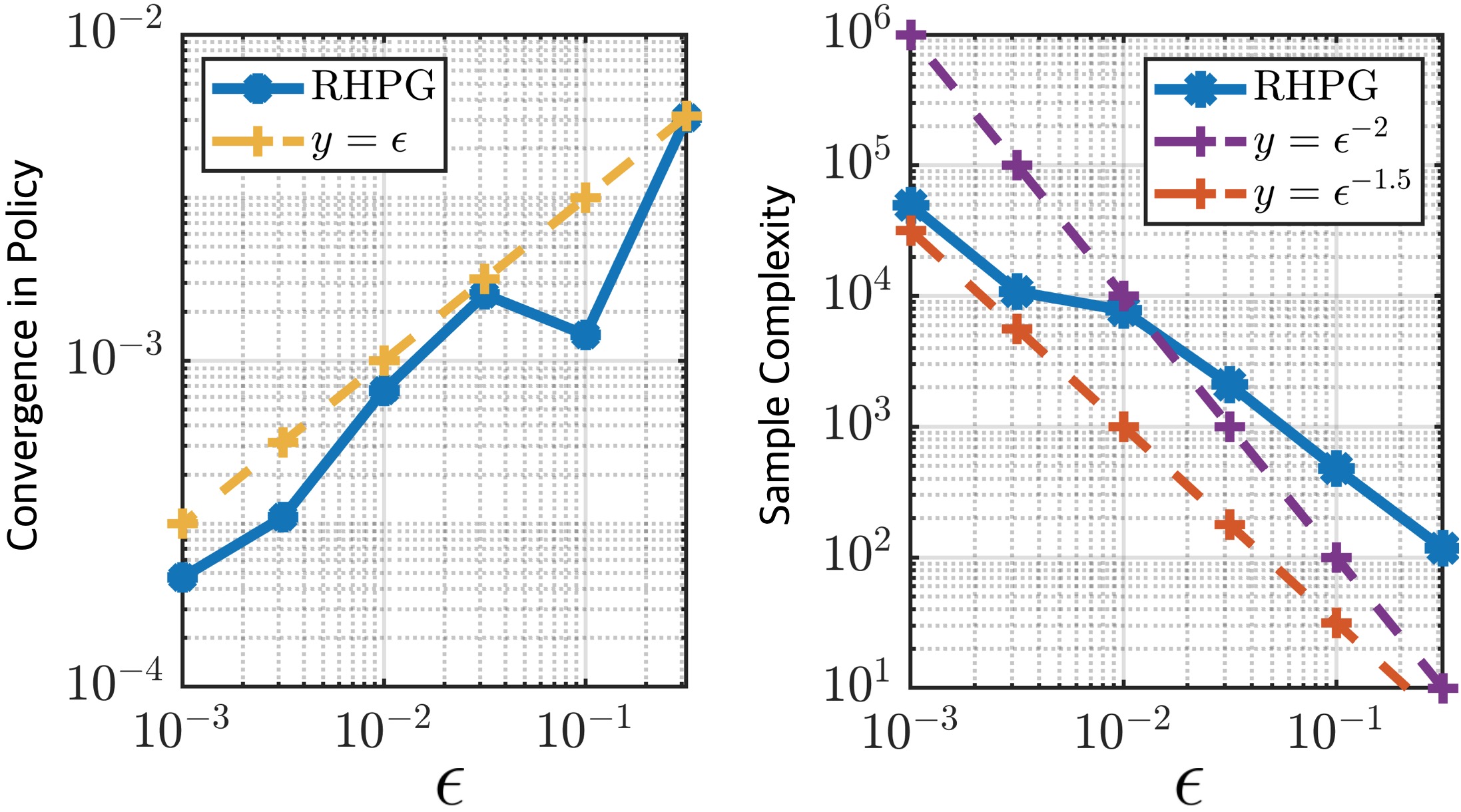}
	\caption{For six different  values of $\epsilon$: \emph{Left}: policy error between the output and $K^*$. \emph{Right}: the total number of calls to the (two-point) zeroth-order oracle.}\label{fig:sim}
\end{figure}

We verify our theories on a scalar linear system studied in \cite{malik2020derivative}, where $A = 5$, $B = 0.33$, $Q = R = 1$, and the unique optimal LQR policy is $K^* = 14.5482$. For the PG method in \cite{fazel2018global, malik2020derivative} to converge in this simple setting, one has to initialize with a policy $K_0$ that satisfies $K_0\in \underline{\cK}:=\{K \mid 12.12 < K \pm r <18.18 \}$, which is necessary to prevent the zeroth-order oracle with a smoothing radius of $r$ from perturbing $K_0$ outside of the stabilizing region $\cK$ during the first iteration of the PG update. In contrast, we initialize the PG updates in Algorithm \ref{alg:DP} with a zero policy $K_h = 0$, set $Q_N = 3$, and choose $N = \texttt{ceil}(\log(\epsilon^{-1}))$ according to \eqref{eqn:N0}. Furthermore, we choose $r_h = \sqrt{\epsilon}$, select a constant stepsize in each iteration of the RHPG algorithm, and run the algorithm to solve the LQR problem under six different $\epsilon$, namely $\epsilon \in \{10^{-3}, 3.16\times 10^{-3}, 10^{-2}, 3.16\times 10^{-2}, 10^{-1}, 3.16\times 10^{-1}\}$. We apply the zeroth-order PG update in solving every subproblem to $\big\|\tilde{K}_h - \tilde{K}^*_h\big\| \leq \epsilon$. As shown in Figure \ref{fig:sim}, the empirical observation of the iteration complexity of RHPG (right) for the convergence in policy (left) is around $\cO(\epsilon^{-2})$ under varying $\epsilon$, which corroborates our theoretical findings.

\section{Discussions}\label{sec:discussion}
In this section, we provide two discussion paragraphs. The first paragraph discusses the tradeoffs in selecting the problem horizon $N$. The second paragraph covers extensions of our results to the stochastic LQR setting that incorporates a zero-mean, independent stochastic disturbance in the dynamics and, additionally, another setting where the initial state $x_0$ could be arbitrary.

\subsection{Tradeoffs in Selecting the Problem Horizon $N$}\label{sec:tradeoff}
There exist two potential usages of the RHPG algorithm in practice. In most real-world scenarios, the user wants to address the finite-horizon LQR problem with \emph{time-varying} system parameters and does not necessarily care whether the resulting control policy is close to the infinite-horizon solution or not. Model-predictive or receding-horizon control is the most prevalent choice when detailed model information is available. The RHPG solution, on the other hand, extends receding-horizon control to the setting with an unknown model, and the user could specify a desired horizon $N$. 

In more theory-oriented cases, the user wants to address the infinite-horizon LQR problem with time-invariant system parameters. Then the user should choose the problem horizon $N$ carefully, balancing between i) reducing the finite-to-infinite horizon error by increasing $N$, and ii) potentially creating more computational error in the approximate dynamic programming when $N$ is large. The RHPG framework performs these two tasks sequentially: first dealing with the finite-to-infinite horizon error and fixing an $N$, and then addressing the finite-$N$-horizon problem step-by-step in time. The main rationale comes from the different rates of $N$ contributing to the two errors. In the first step (cf., Theorem \ref{lemma:finite_approximation}), the convergence of the finite-horizon solution toward the infinite-horizon solution is exponential, meaning that it suffices to choose a problem horizon $N = \cO(\log(\epsilon^{-1}))$ for an accurate approximation (the dependence to other system parameters are also logarithmic). After deciding on an $N$ based on the desired $\epsilon$, the user solves the finite-$N$-horizon problem iteratively, where the subproblem in each iteration is strongly convex and smooth. In the sample complexity analysis, we have considered the choice of $N$ since one would require a higher solution accuracy in each iteration if $N$ is large (cf., \S\ref{proof:LQR_DP}). Moreover, $N$ also appears in the probability term of $1-N\delta$ due to the utilization of Boole's inequality. In solving each sub-problem, the complexity in the failure probability $\delta$ is $\cO(\log(\delta^{-1}))$. Replacing $\delta$ with $\delta/N$ yields the rate of $\cO(\log(N/\delta)) = \cO(\log(\delta^{-1})) + \cO(\log\log(\epsilon^{-1}))$. Hence, with a fixed $N$ chosen according to Theorem \ref{lemma:finite_approximation}, the impact of $N$ on the total complexity is an additional $\cO(\log\log(\epsilon^{-1}))$ factor and is therefore minor. In summary, the user in the infinite-horizon setting should prioritize the choice of $N$ toward reducing the finite-to-infinite horizon error according to Theorem \ref{lemma:finite_approximation} and then solve the finite-$N$-horizon problem.

\subsection{Extensions to Arbitrary $x_0$ and Stochastic LQR}\label{sec:extension}
We discuss here extensions of the RHPG framework to i) the setting where $x_0\in \RR^n$ is an arbitrary (deterministic) vector that is unknown to the designer, and ii) stochastic LQR with $x_0\sim\cD$ and the system dynamics being
\begin{align*}
	x_{t+1} = Ax_t + Bu_t + w_t, \quad w_t \sim \cN(0, W), \quad W>0.
\end{align*} 
Note that in the stochastic LQR setting, the objective function \eqref{eqn:cost} should be replaced with the time-average cost
\begin{align*}
		 J_{\infty} := \limsup_{N\to\infty}\frac{1}{N}~\EE_{x_0, w_t}\hspace{-0.1em}\left[\sum_{t=0}^{N-1} \big( x_t^{\top}Qx_t + u^{\top}_tRu_t \big)\right].
\end{align*}
In both settings and under the stabilizability and detectability assumptions, a unique stationary state-feedback control policy exists and is identical to \eqref{eqn:infinite_lqr_gain}, which also stabilizes the closed-loop system. The difference is that in setting i), the LQR problem is deterministic, which allows implementing RHPG with a two-point zeroth-order oracle (cf., Algorithm \ref{alg:DP}). In contrast, setting ii) involves additive noises in the system dynamics, which necessitates using a one-point zeroth-order oracle, and thus, the gradient sampling will be noisier.

We note that the procedure of RHPG is identical whether solving deterministic or stochastic problems and with either finite or infinite problem horizons. The only difference is in the choice of the inner-loop oracles (two-point v.s. one-point). In both settings, the receding-horizon parametrization, the convergence of the Riccati equations in Theorem \ref{lemma:finite_approximation}, the analysis of dynamic programming in Theorem \ref{theorem:LQR_DP}, and the quadratic optimization landscape in each subproblem are the same as in the presentation in the main body of this paper. Simply combining Theorems \ref{lemma:finite_approximation}-\ref{theorem:LQR_DP} with the corresponding inner-loop convergence result (for quadratic minimization and as a replacement of Proposition 3.3) yields the overall complexity in these two extended settings.

\section{Conclusion}
We have revisited discrete-time LQR from the perspective of RHPG and provided a fine-grained sample complexity analysis for RHPG to learn a control policy that is stabilizing and $\epsilon$-close to the optimal LQR policy. Our result demonstrates the potential of RHPG in addressing various tasks in linear control and estimation with streamlined analyses.

\bibliographystyle{unsrt}
\bibliography{main}

\renewcommand{\theequation}{\thesubsection.\arabic{equation}}

\appendix
\subsection{Proof of Theorem \ref{lemma:finite_approximation}}\label{proof:finite}
This proof is dual to the proof of Theorem 3.1 in \cite{zhang2023learning}. We first identify one-to-one correspondences between system parameters in LQR and those in Kalman filtering \cite{zhang2023learning}:
\begin{center}
	\begin{tabular}{|l|l|l|l|l|l|}
\cline{1-6}
LQR: & $A$ & $B$ & $Q$ & $R$ & $Q_N$ \\ 
     & \hspace{0.2em}$\updownarrow$ & \hspace{0.2em}$\updownarrow$ & \hspace{0.2em}$\updownarrow$ & \hspace{0.2em}$\updownarrow$ & \hspace{0.2em}$\updownarrow$\\ 
KF \cite{zhang2023learning}:  & $A^{\top}$ & $C^{\top}$ & $W$ & $V$ & $X_0$\\ \cline{1-6}             
\end{tabular}
\end{center}
We also identify direct correspondences between our $P_t$, $P^*$, $K_t$, and $K^*$ and \cite{zhang2023learning}'s $\Sigma_{N-t}$, $\Sigma^*$, $L_{N-1-t}$, and $L^*$, respectively. Then, letting 
\begin{align*}
	&\tilde{P}_t := P^*_t- P^*\hspace{-0.2em}, \quad  \tilde{R} := R + B^{\top}P^*B, \\
	 &\overline{A} := A - B\tilde{R}^{-1}B^{\top}P^*A,
\end{align*}
and following equations (A.2)-(A.3) of \cite{zhang2023learning}, we have
\begin{align}
	\tilde{P}_{t} &= \overline{A}^{\top}\tilde{P}_{t+1}\overline{A} - \overline{A}^{\top}\hspace{-0.1em}\tilde{P}_{t+1}B(\tilde{R} + B^{\hspace{-0.1em}\top}\hspace{-0.1em}\tilde{P}_{t+1}B)^{\hspace{-0.1em}-1}\hspace{-0.1em}B^{\hspace{-0.1em}\top}\hspace{-0.1em}\tilde{P}_{t+1}\overline{A} \nonumber\\
	&= \overline{A}^{\top}\tilde{P}_{t+1}^{1/2}\big[I + \tilde{P}^{1/2}_{t+1}B\tilde{R}^{-1}B^{\top}\tilde{P}_{t+1}^{1/2}\big]^{-1}\tilde{P}_{t+1}^{1/2}\overline{A} \nonumber\\
	&\leq \hspace{-0.1em}\big[1\hspace{-0.15em}+\hspace{-0.15em}\lambda_{\min}(\tilde{P}^{1/2}_{t+1}B\tilde{R}^{-1}B^{\top}\tilde{P}_{t+1}^{1/2})\big]^{-1}\hspace{-0.1em}\overline{A}^{\top}\tilde{P}_{t+1}\overline{A} \nonumber\\
	&=:\hspace{-0.1em} \mu_{t} \overline{A}^{\top}\hspace{-0.2em}\tilde{P}_{t+1}\overline{A},\label{eqn:RDE_conv_step2}
\end{align}
where $\tilde{P}_{t+1}^{1/2}$ denotes the unique positive semi-definite (psd) square root of the psd matrix $\tilde{P}_{t+1}$, $0 < \mu_t \leq 1$ for all $t$, and $\overline{A}$ satisfies $\rho(\overline{A}) < 1$. We now use $\|\cdot\|_*$ to represent the $P^*$-induced matrix norm and invoke Theorem 14.4.1 of \cite{hassibi1999indefinite}, where our $\tilde{P}_t$, $\overline{A}^{\top}$ and $P^*$ correspond to $P_i - P^*$, $F_p$ and $W$ in \cite{hassibi1999indefinite}, respectively. By Theorem 14.4.1 of \cite{hassibi1999indefinite} and \eqref{eqn:RDE_conv_step2}, we obtain $\|\overline{A}\|_* < 1$ and given that $\mu_t \leq 1$, 
\begin{align*}
	\|\tilde{P}_{t}\|_* \leq \|\overline{A}\|^2_* \cdot \|\tilde{P}_{t+1}\|_*.
\end{align*}
Therefore, the convergence is exponential such that $\|\tilde{P}_t\|_* \leq \|\overline{A}\|_*^{2(N-t)}\cdot \|\tilde{P}_{N}\|_*$. As a result, the convergence of $\tilde{P}_t$ to $0$ in spectral norm can be characterized as
\begin{align*}
	\|\tilde{P}_t\| \leq \kappa_{P^*}\cdot \|\tilde{P}_t\|_* \leq \kappa_{P^*}\cdot\|\overline{A}\|_*^{2(N-t)}\cdot \|\tilde{P}_{N}\|_*,
\end{align*}
where we have used $\kappa_X$ to denote the condition number of $X$. That is, to ensure $\|\tilde{P}_1\| \leq \epsilon$, it suffices to require
\begin{align}\label{eqn:required_time}
	N \geq \frac{1}{2}\cdot \frac{\log\big(\frac{\|\tilde{P}_N\|_*\cdot \kappa_{P^*}}{\epsilon}\big)}{\log\big(\frac{1}{\|\overline{A}\|_*}\big)} + 1.
\end{align}
Lastly, we show that the (monotonic) convergence of $K^*_t$ to $K^*$ follows from the convergence of $P^*_t$ to $P^*$. Similar to (A.5) of \cite{zhang2023learning}, this can be verified through:
\begin{align}
	K^*_t - K^* &= (R+B^{\top}P^*_{t+1}B)^{-1}B^{\top}P^*_{t+1}A \nonumber\\
	&\hspace{1em}- (R+B^{\top}P^*B)^{-1}B^{\top}P^*A \nonumber\\
	&\hspace{-3.7em}=\big[(R+B^{\top}P^*_{t+1}B)^{-1}-(R+B^{\top}P^*B)^{-1}\big]B^{\top}P^*A \nonumber\\
	&\hspace{-2.5em} + (R+B^{\top}P^*_{t+1}B)^{-1}B^{\top}(P^*_{t+1}-P^*)A \nonumber\\
	&\hspace{-3.7em}= (R+B^{\top}P^*_{t+1}B)^{-1}B^{\top}(P^*-P^*_{t+1})BK^* \nonumber\\
	&\hspace{-2.5em}-(R+B^{\top}P^*_{t+1}B)^{-1}B^{\top}(P^*-P^*_{t+1})A\nonumber\\
	&\hspace{-3.7em}=(R+B^{\top}P^*_{t+1}B)^{-1}B^{\top}(P^*-P^*_{t+1})(BK^*-A) \label{eqn:kdiff}
\end{align}
Hence, we have $\|K^*_t - K^*\| \leq \frac{\|\overline{A}\|\cdot \|B\|}{\lambda_{\min}(R)}\cdot \|P^*_{t+1} - P^*\|$ and
\begin{align*}
	\|K^*_0 - K^*\| \leq \frac{\|\overline{A}\|\cdot \|B\|}{\lambda_{\min}(R)}\cdot \|\tilde{P}_1\|.
\end{align*}
Substituting $\epsilon$ in \eqref{eqn:required_time} with $\frac{\epsilon\cdot\lambda_{\min}(R)}{\|\overline{A}\|\cdot\|B\|}$ completes the proof.

\subsection{Proof of Theorem \ref{theorem:LQR_DP}}\label{proof:LQR_DP}
This proof is dual to the proof of Theorem 3.3 in \cite{zhang2023learning}. First, according to Theorem \ref{lemma:finite_approximation}, we select
\begin{align}\label{eqn:N_choice}
	N = \frac{1}{2}\cdot \frac{\log\big(\frac{2\|Q_N-P^*\|_*\cdot\kappa_{P^*}\cdot \|A_K^*\|\cdot\|B\|} {\epsilon\cdot\lambda_{\min}(R)}\big)}{\log\big(\frac{1}{\|A_K^*\|_*}\big)} + 1,
\end{align} 
where $A_K^*:=A-BK^*$. This ensures that $K^*_{0}$ is stabilizing and $\|K^*_{0} - K^*\| \leq \epsilon/2$. Then, it remains to show that the output $\tilde{K}_0$ satisfies $\|\tilde{K}_0 - K^*_0\| \leq \epsilon/2$. 

Recall that the RDE \eqref{eqn:RDE} is a backward iteration starting with $P^*_N = Q_N \geq 0$, and can also be represented as: 
\begin{align}
	P^*_{t} &= A^{\top}P^*_{t+1}\big(A - BK^*_t\big) + Q \label{eqn:standard_RDE_appen}\\
	&\hspace{-1em}=(A\hspace{-0.1em}-\hspace{-0.1em}BK^*_t)^{\hspace{-0.1em}\top}\hspace{-0.1em}P^*_{t+1}(A\hspace{-0.1em}-\hspace{-0.1em}BK^*_t) \hspace{-0.1em}+\hspace{-0.1em} (K^*_t)^{\hspace{-0.1em}\top}\hspace{-0.1em}RK^*_t \hspace{-0.1em}+\hspace{-0.1em} Q \label{eqn:lqr_RDE_Lya}. 
\end{align}
 Moreover, for any $K_t$, we introduce the Lyapunov equation:
 \begin{align}
 	\hspace{-0.5em}P_{t} = (A-BK_t)^{\top}P_{t+1}(A-BK_t) + K_t^{\top}RK_t+Q. \label{eqn:lqr_lyapunov}
 \end{align}
 Furthermore, for clarity of the proof, we define/recall:
\small
\begin{align*}
 	&K^*_t\text{: Exact LQR policy at time $t$ defined in \eqref{eqn:finite_lqr_gain}}\\
 	&\tilde{K}_t^*\text{: Optimal policy of the current cost-to-go function,}\\
 	&\hspace{2.1em} \text{absorbing errors in all prior steps} \\
 	&\tilde{K}_t\text{: An approximation of $\tilde{K}_t^*$ obtained by the PG update \eqref{eqn:PG}}\\
 	&\hspace{0.2em} \delta_t:=\tilde{K}_t-\tilde{K}_t^* \text{: Policy optimization error at time $t$} \\
 	&\tilde{P}^*_{t} \text{: Generated by \eqref{eqn:lqr_RDE_Lya} with $K^*_t = \tilde{K}^*_t$ and $P^*_{t+1} = \tilde{P}_{t+1}$.}
 \end{align*}
 \normalsize
 
 We argue that $\|\tilde{K}_{0} - K^*_{0}\| \leq \epsilon/2$ can be achieved by carefully controlling $\delta_t$ for all $t$. At $t=0$, it holds that
 \begin{align*}
 	\|\tilde{K}_{0} - K^*_{0}\| \leq \|\tilde{K}^*_{0} - K^*_{0}\| + \|\delta_{0}\|,
 \end{align*}
 where substituting $K^*_t$ and $K^*$ in \eqref{eqn:kdiff}, respectively, with $\tilde{K}^*_0$ and $K^*_0$ leads to
 \begin{align*}
 	\tilde{K}^*_{0} - K^*_{0} = (R+B^{\top}\tilde{P}_{1}B)^{-1}B^{\top}(P^*_1-\tilde{P}_{1})(BK^*_0-A).
 \end{align*}
 Hence, the error size $\|\tilde{K}^*_{0} - K^*_{0}\|$ could be bounded by 
 \begin{align}
	\|\tilde{K}^*_{0} - K^*_{0}\| \leq \frac{\|A-BK^*_0\|\cdot\|B\|}{\lambda_{\min}(R)}\cdot\|P^*_1-\tilde{P}_{1}\| \label{eqn:laststep_req}.
\end{align}
Define the helper constants
\begin{align*}
	C_1:= \frac{\varphi\cdot\|B\|}{\lambda_{\min}(R)} > 0, \ \varphi := \max_{t\in\{0, \cdots, N-1\}}\|A-BK^*_t\|.
\end{align*}
Next, we require $\|\delta_{0}\| \leq \epsilon/4$ and $\|\tilde{K}^*_{0} - K^*_{0}\| \leq \epsilon/4$ to fulfill $\|\tilde{K}_{0} - K^*_{0}\| \leq \epsilon/2$. We select a fixed scalar $a > 0$ that is independent of system parameters and $\epsilon$, and additionally require $\|P^*_1-\tilde{P}_{1}\| \leq a$ to upper-bound the pd solutions of \eqref{eqn:lqr_lyapunov}. Then, by \eqref{eqn:laststep_req}, in order to fulfill $\|\tilde{K}^*_{0} - K^*_{0}\| \leq \epsilon/4$, it suffices to require
\begin{align}\label{eqn:preq1}
\|P^*_1-\tilde{P}_{1}\| \leq \min\bigg\{a, \frac{\epsilon}{4 C_1}\bigg\}.
\end{align}
Subsequently, by \eqref{eqn:lqr_lyapunov}, we have
\begin{align}\label{eqn:plast}
	P^*_1-\tilde{P}_{1} = (P^*_1 - \tilde{P}^*_{1}) + (\tilde{P}^*_{1} - \tilde{P}_{1}).
\end{align}
The first difference term on the RHS of \eqref{eqn:plast} is
\begin{align}
P^*_1 - \tilde{P}^*_{1} &= A^{\top}P^*_{2}\big(A - BK^*_1\big) - A^{\top}\tilde{P}_{2}\big(A - B\tilde{K}^*_1\big) \nonumber\\
&\hspace{-4.4em}= A^{\top}(P_2^* - \tilde{P}_2)(A-BK^*_1) + A^{\top}\tilde{P}_2B(\tilde{K}_1^* - K^*_1).\label{eqn:plast2}\\
&\hspace{-4.4em}=A^{\top}(P_2^* - \tilde{P}_2)(A-BK^*_1) \nonumber\\
&\hspace{-3.4em}- A^{\top}\tilde{P}_2B(R+B^{\top}\tilde{P}_{2}B)^{-1}B^{\top}(P^*_2-\tilde{P}_{2})(A-BK^*_1) \label{eqn:apply_k_diff}\\
&\hspace{-4.4em}=A^{\top}[I-\tilde{P}_2B(R+B^{\top}\tilde{P}_{2}B)^{-1}B^{\top}](P_2^* - \tilde{P}_2)(A-BK^*_1) \nonumber \\
&\hspace{-4.4em}=A^{\top}(I+\tilde{P}_2BR^{-1}B^{\top})^{-1}(P_2^* - \tilde{P}_2)(A-BK^*_1), \label{eqn:dp_inversion_lemma}
\end{align}
where applying  \eqref{eqn:laststep_req} in deriving \eqref{eqn:apply_k_diff} and \eqref{eqn:dp_inversion_lemma} is due to the matrix inversion lemma.
Moreover, the second term on the RHS of \eqref{eqn:plast} is 
\begin{align}
	&\tilde{P}^*_{1} - \tilde{P}_{1} = (A-B\tilde{K}^*_1)^{\top}\tilde{P}_{2}(A-B\tilde{K}^*_1) + (\tilde{K}^*_1)^{\top}R\tilde{K}^*_1 \nonumber\\
	&\hspace{5em} - (A-B\tilde{K}_1)^{\top}\tilde{P}_{2}(A-B\tilde{K}_1) - (\tilde{K}_1)^{\top}R\tilde{K}_1 \nonumber\\
	&= - (\tilde{K}^*_1)^{\hspace{-0.1em}\top}\hspace{-0.1em}B^{\hspace{-0.1em}\top}\hspace{-0.1em}\tilde{P}_2A \hspace{-0.1em}-\hspace{-0.1em} A^{\hspace{-0.1em}\top}\hspace{-0.1em}\tilde{P}_2B\tilde{K}_1^* \hspace{-0.1em}+\hspace{-0.1em} (\tilde{K}^*_1)^{\hspace{-0.1em}\top}\hspace{-0.1em}(R+B^{\hspace{-0.1em}\top}\tilde{P}_2B)\tilde{K}^*_1   \nonumber\\
	&\hspace{1em} + \tilde{K}_1^{\top}B^{\top}\tilde{P}_2A + A^{\top}\tilde{P}_2B\tilde{K}_1 - \tilde{K}_1^{\top}(R+B^{\top}\tilde{P}_2B)\tilde{K}_1 \nonumber\\
	&= \big[(R+B^{\top}\tilde{P}_2B)^{-1}B^{\top}\tilde{P}_2A - \tilde{K}^*_1\big]^{\top}(R+B^{\top}\tilde{P}_2B)\cdot \nonumber \\
	&\hspace{1.5em} \big[(R+B^{\top}\tilde{P}_2B)^{-1}B^{\top}\tilde{P}_2A - \tilde{K}^*_1 \big] \nonumber\\
	&\hspace{1em}- \big[(R+B^{\top}\tilde{P}_2B)^{-1}B^{\top}\tilde{P}_2A - \tilde{K}_1\big]^{\top}(R+B^{\top}\tilde{P}_2B)\cdot \nonumber\\
	&\hspace{1.5em}\big[(R+B^{\top}\tilde{P}_2B)^{-1}B^{\top}\tilde{P}_2A - \tilde{K}_1 \big] \label{eqn:complete_square}\\
	&= - \delta_1^{\top}(R+B^{\top}\tilde{P}_2B)\delta_1, \label{eqn:pdiff2}
\end{align}
where \eqref{eqn:complete_square} follows from completion of squares. Thus, combining  \eqref{eqn:plast}, \eqref{eqn:plast2}, and \eqref{eqn:pdiff2} yields
\begin{align}
	&\hspace{1.3em}\|P^*_1-\tilde{P}_{1}\| \nonumber\\
	&\leq \|P^*_{2} -\tilde{P}_{2}\|\cdot\varphi\|A\| \|(I+\tilde{P}_2BR^{-1}B^{\top})^{-1}\| \nonumber\\
	&\hspace{1em}+ \|\delta_{1}\|^2\|R+B^{\top}\tilde{P}_2B\| \nonumber \\
	&\leq \varphi\|A\| \cdot\|P^*_{2} - \tilde{P}_{2}\| + \|\delta_{1}\|^2\|R+B^{\top}\tilde{P}_2B\| \label{eqn:p_contraction},
\end{align}
where \eqref{eqn:p_contraction} is due to that $\tilde{P}_2BR^{-1}B^{\top} \geq 0$ and thus $\|(I+\tilde{P}_2BR^{-1}B^{\top})^{-1}\| \leq 1$. Now, we require
\begin{align}
	\|P^*_{2} - \tilde{P}_{2}\| &\leq \min\bigg\{a, \frac{a}{C_2}, \frac{\epsilon}{4 C_1C_2}\cdot\bigg\} \label{eqn:preq2}\\
	\|\delta_{1}\| &\leq  \min\bigg\{\sqrt{\frac{a}{C_3}}, \frac{1}{2}\sqrt{\frac{\epsilon}{C_1C_3}}\bigg\}\label{eqn:K_req1},
\end{align}
where $C_2$ and $C_3$ are positive constants defined as\footnote{As the scalar $a>0$ increases, the constant $C_3$ grows correspondingly.}
\begin{align*}
	&C_2 := 2\varphi\|A\| >0, \ C_3 := 2\|R+B^{\top}(P_{\max} + aI)B\| > 0 \\ 
	&P_{\max} := \max_{t\in\{0, \cdots, N-1\}}\{P^*_t\}.
\end{align*}
Then, conditions \eqref{eqn:preq2} and \eqref{eqn:K_req1} are sufficient for \eqref{eqn:preq1} (and thus for $\|\tilde{K}_{0} - K^*_{0}\| \leq \epsilon/2$) to hold. Subsequently, we can propagate the required accuracies in \eqref{eqn:preq2} and \eqref{eqn:K_req1} forward in time. Specifically,   we iteratively apply the arguments in \eqref{eqn:p_contraction} (i.e., by plugging quantities with subscript $t$ into the LHS of \eqref{eqn:p_contraction} and plugging quantities with subscript $t+1$ into the RHS of \eqref{eqn:p_contraction}) to obtain  the result that if at all $t \in \{2, \cdots, N-1\}$, we require
\begin{align}
	&\|P^*_{t} - \tilde{P}_{t}\| \leq \min\bigg\{a, \frac{a}{C_2^{t-1}}, \frac{\epsilon}{4 C_1C_2^{t-1}}\bigg\} \label{eqn:preq_allt}\\
	&\hspace{-0.5em}\|\delta_{t}\| \leq  \min\bigg\{\sqrt{\frac{a}{C_3}}, \sqrt{\frac{a}{C_2^{t-2}C_3}}, \frac{1}{2}\sqrt{\frac{\epsilon}{C_1C_2^{t-2}C_3}}\bigg\} \nonumber,
\end{align}
then \eqref{eqn:preq2} holds true and therefore \eqref{eqn:preq1} is satisfied.

We now compute the required accuracy for $\delta_{N-1}$. Note that $P^*_{N-1} = \tilde{P}^*_{N-1}$ since no prior computational errors happened at $t=N$. By \eqref{eqn:p_contraction}, the distance between $P^*_{N-1}$ and $\tilde{P}_{N-1}$ can be bounded as 
\begin{align*}
	\|P^*_{N-1} - \tilde{P}_{N-1}\| = \|\tilde{P}^*_{N-1} - \tilde{P}_{N-1}\| \leq \|\delta_{N-1}\|^2\cdot C_3.
\end{align*}
To fulfill the requirement \eqref{eqn:preq_allt} for $t=N-1$, which is
\begin{align*}
	\|P^*_{N-1} - \tilde{P}_{N-1}\| \leq \min\bigg\{a, \frac{a}{C_2^{N-2}}, \frac{\epsilon}{4 C_1C_2^{N-2}}\bigg\},
\end{align*}
it suffices to let
\small
\begin{align}
	\hspace{-0.4em}\|\delta_{N-1}\| \leq \min \hspace{-0.1em}\bigg\{\sqrt{\frac{a}{C_3}}, \sqrt{\frac{a}{C_2^{N-2}C_3}}, \frac{1}{2}\sqrt{\frac{\epsilon}{C_1C_2^{N-2}C_3}}\bigg\}\hspace{-0.1em}. \label{eqn:delta0_req}
\end{align}
\normalsize

Finally, we analyze the worst-case complexity of RHPG by computing, at the most stringent case, the required size of $\|\delta_t\|$. When $C_2 \leq 1$, the most stringent dependence of $\|\delta_t\|$ on $\epsilon$ happens at $t=0$, which is of the order $\cO(\epsilon)$, and the dependences on system parameters are $\cO(1)$. We then analyze the case where $C_2 > 1$, where the requirement on $\|\delta_0\|$ is still $\cO(\epsilon)$. Note that in this case, $\|\delta_{N-1}\| \leq \|\delta_t\|$ for all $t \in \{1, \cdots, N-1\}$ and by \eqref{eqn:delta0_req}:
\begin{align}\label{eqn:delta0}
	\|\delta_{N-1}\| \sim \cO\Big(\sqrt{\frac{\epsilon}{C_1C_2^{N-2}C_3}}\Big).
\end{align}
Since we require $N$ to satisfy \eqref{eqn:N_choice}, the dependence of $\|\delta_{N-1}\|$ on $\epsilon$ in \eqref{eqn:delta0} becomes $\|\delta_{N-1}\| \sim \cO(\epsilon^{\frac{3}{4}})$ with additional polynomial dependences on system parameters, but one can observe that the dependence on $\epsilon$ is still milder than the requirement for $\|\delta_{0}\|$. Therefore, it suffices to require error bound for all $t$ to be $\|\delta_t\| \sim \cO(\epsilon)\cO(\texttt{poly}(\text{system parameters)})$ to reach the $\epsilon$-neighborhood of the infinite-horizon LQR policy. Lastly, for $\tilde{K}_{0}$ to be stabilizing, it suffices to select a sufficiently small $\epsilon$ such that the $\epsilon$-ball centered at the infinite-horizon LQR policy $K^*$ lies entirely in the set of stabilizing policies. A crude bound that satisfies this requirement is
\begin{align*}
	\epsilon < \frac{1 - \|A-BK^*\|_*}{\|B\|} \Longrightarrow \|A-B\tilde{K}_0\|_* < 1.
\end{align*}
This completes the proof.

\subsection{Proof of Proposition \ref{prop:sample}}\label{proof:sample}
Recall that for all $h$, the objective function $ J_h$ is $L_h$-smooth and $\alpha_h$-strongly-convex. Define $\varsigma_h := \frac{\epsilon^2 \alpha_h}{2}$ and $\varsigma:=\min_h \varsigma_h > 0$. We argue that if with a probability of at least $1-\delta$, it holds that
	\begin{align}\label{eqn:cost_convergence}
		 J_h(K_{h, T_h}) -  J_h(\tilde{K}_{h}^*) \leq \varsigma,
	\end{align}
	then $\|K_{h, T_h} - \tilde{K}^*_h\| \leq \epsilon$ also holds with a probability of at least $1-\delta$ and the proof of Proposition \ref{prop:sample} is complete. This is due to the $\alpha_h$-strong convexity and $\nabla  J_h(\tilde{K}^*_h) = 0$. Thus,
	\begin{align*}
		&\hspace{1em} J_h(K_{h, T_h}) -  J_h(\tilde{K}_h^*) \\
		&\geq  \nabla  J_h(\tilde{K}^*_h)^{\top} (K_{h, T_h} - \tilde{K}_h^*) + \frac{\alpha_h}{2} \|K_{h, T_h} - \tilde{K}^*_h\|^2_F \\
		&\Longrightarrow  \|K_{h, T_h} - \tilde{K}^*_h\|^2_F \leq \frac{2}{\alpha_h} \big[ J_h(K_{h, T_h}) -  J_h(\tilde{K}_h^*)\big] \leq \epsilon. 
	\end{align*}
As a result, we will focus on proving \eqref{eqn:cost_convergence} with a high probability of at least $1-\delta$. First, define the cost difference $\Delta_t:= J_h(K_{h, t}) -  J_h(\tilde{K}_{h}^*)$ and the stopping time $\tau:=\min\{t \mid \Delta_t > 10\delta^{-1}\Delta_0\}$. Let $\EE^t[\cdot]$ denote the expectation conditioned on all the randomness up to $t$. Then, we state the following helper lemma and defer its proof to \S\ref{proof:zeroth_descent}.

\begin{lemma}\label{lemma:zeroth_descent}
For all $h$, choose the parameters of Algorithm \ref{alg:DP} according to 
\begin{align*}
	\eta_h \leq \frac{1}{2L_h}, \quad r_h \leq \min\Big\{\frac{\alpha_h}{4L_h}\sqrt{\frac{\varsigma\delta}{10}}, \frac{1}{2L_h}\sqrt{\frac{\alpha_h\varsigma\delta}{5}}\Big\}.
\end{align*}
Then, for all $t$, it holds that
	\begin{align*}
		\EE^t[\Delta_{t+1}] \leq \Big(1-\frac{\eta_h\alpha_h}{4}\Big)\Delta_t + \frac{L_h\eta_h^2}{2}G_2 + \eta_h\alpha_h\frac{\varsigma\delta}{20},
	\end{align*}
	where $\alpha_h$ and $L_h$ are the strong convexity and smoothness constants of $ J_h$, respectively, and $G_2$ is a uniform constant to be introduced shortly. 
\end{lemma} 

Following the proof of Theorem 8 in \cite{malik2020derivative, furieri2019learning}, we first consider the case of $\tau > T_h$. In this case, we can bound $\EE^t[\Delta_{t+1}]$ using Lemma \ref{lemma:zeroth_descent} directly. When $\tau \leq T_h$, it implies that $\EE^t[\Delta_{t+1}]1_{\tau>t} = 0$. We require $\eta_h \leq \frac{\varsigma\delta\alpha_h}{40L_hG_2}$ and show
\begin{align*}
	\EE^t[\Delta_{t+1}]1_{\tau>t+1} &\leq \Big(1-\frac{\eta_h\alpha_h}{4}\Big)^{t+1}\Delta_0 \\
	&+ \Big(\frac{L_h\eta_h^2}{2}G_2 \hspace{-0.1em}+\hspace{-0.1em} \eta_h\alpha_h\frac{\varsigma\delta}{20}\Big)\hspace{-0.1em}\sum_{i=0}^{t}\hspace{-0.1em}\Big(1\hspace{-0.1em}-\hspace{-0.1em}\frac{\eta_h\alpha_h}{4}\Big)^{\hspace{-0.1em}i} \\
	&\leq \Big(1-\frac{\eta_h\alpha_h}{4}\Big)^{t+1}\Delta_0 + \frac{\varsigma\delta}{4}
\end{align*}
Setting $t+1 =T_h$, it suffices to let $T_h = \frac{4}{\eta_h\alpha_h}\log(\frac{4\Delta_0}{\delta\varsigma})$ to ensure that
\begin{align*}
	\EE[\Delta_{T_h}1_{\tau>T_h}] \leq \Big(1-\frac{\eta_h\alpha_h}{4}\Big)^{T_h}\Delta_0 + \frac{\varsigma\delta}{4} \leq \frac{\varsigma\delta}{2}.
\end{align*}
Next, we prove that the event $\tau \leq T_h$ has a probability smaller than $\frac{\delta}{2}$. For all $t$, we define the stopping process as 
\begin{align*}
	Y_t:=\Delta_{\min\{\tau, t\}} + (T_h - t)\Big(\frac{L_h\eta_h^2}{2}G_2 +\eta_h\alpha_h\frac{\varsigma\delta}{20}\Big),
\end{align*}
By Eq. (20)-(21) of \cite{malik2020derivative}, $Y_t$ is a super-martingale. Applying Doob's maximal inequality yields
\begin{align*}
	&P\Big(\max_{t=1, \cdots, T_h}Y_t \geq \frac{10\Delta_0}{\delta}\Big) \leq \frac{\delta\EE[Y_0]}{10\Delta_0}\\
	&=\frac{\delta}{10\Delta_0}\Big(\Delta_0 + T_h\big(\frac{L_h\eta_h^2}{2}G_2 + \eta_h\alpha_h\frac{\varsigma\delta}{20}\big)\Big)\\
	&\leq \frac{\delta}{10\Delta_0}\Big(\Delta_0 + \log\Big(\frac{4\Delta_0}{\varsigma\delta}\Big)\frac{\varsigma\delta}{20} + \log\Big(\frac{4\delta_0}{\varsigma\delta}\Big)\frac{\varsigma\delta}{5}\Big)
\end{align*}
Imposing the condition that $\varsigma\log(\frac{4\Delta_0}{\varsigma\delta}) \leq 16\delta^{-1}\Delta_0$, we can prove that $P\big(\max_{t=1, \cdots, T_h}Y_t \geq \frac{10\Delta_0}{\delta}\big) \leq \frac{\delta\EE[Y_0]}{10\Delta_0} \leq \frac{\delta}{2}$. We can now conclude that $\EE[\Delta_{T_h}1_{\tau>T_h}] \leq \frac{\delta\varsigma}{2}$ and the event $\tau$ occurs after $T_h$ with probability at least $1-\frac{\delta}{2}$. As a result, 
\begin{align*}
	P(\Delta_{T_h} \geq \varsigma) &\leq P(\Delta_{T_h} 1_{\tau >T_h} \geq \varsigma) + P(1_{\tau \leq T_h}) \\
	&\leq \frac{\EE[\Delta_{T_h} 1_{\tau > T_h}]}{\varsigma} + P(1_{\tau \leq T_h})\leq \frac{\delta}{2} + \frac{\delta}{2} = \delta,
\end{align*}
where we have used Markov's inequality. This verifies \eqref{eqn:cost_convergence}, and thus $\|K_{h, T_h} - \tilde{K}^*_h\| \leq \epsilon$ is satisfied with a probability of at least $1-\delta$. Lastly, we analyze the constant $G_2$ following Corollary 10 of \cite{malik2020derivative} and \cite{shamir2017optimal}, where
\begin{align*}
	&G_2 = \hspace{-0.3em}\sup_{K_h \in \Phi_h} \hspace{-0.2em}\EE\left[\left\|\frac{mn}{2r_h}[J(K_h+r_hU) - J(K_h-r_hU)]U\right\|_F^2\right]\\
	&\Phi_h:=\{K_h\mid  J(K_h) -  J(\tilde{K}^*_h) \leq 10\delta^{-1}\Delta_0\}.
\end{align*}
By Corollary 10 of \cite{malik2020derivative}, it holds almost surely that $G_2 \leq (mn)\lambda^2$, where $\lambda:=\max_h\lambda_h$ and $\lambda_h$ is the Lipschitz continuity constant of $ J_h$ taken over the compact domain $\Phi_h$. In summary, for $\|K_{h, T_h} - \tilde{K}^*_h\| \leq \epsilon$ to hold with a probability of at least $1-\delta$, we need to choose the (constant) algorithmic parameters according to $\eta_h \sim \cO(\epsilon^2)$, and $r_h \sim \cO(\epsilon)$. Then, the iteration complexity for the convergence of the zeroth-order PG method is  
$T_h \sim \cO(\frac{1}{\epsilon^2}\log(\frac{1}{\delta\epsilon^2}))$.

\subsection{Proof of Lemma \ref{lemma:zeroth_descent}}\label{proof:zeroth_descent}
The proof of this lemma mostly follows the steps in Section 4.1.1 of \cite{malik2020derivative}. First, define the smoothed version of $ J_h$ as $ J_h^{r_h}(K_h):= \EE[ J_h(K_h + r_hU)]$, where the expectation is taken over $U$ that is uniformly drawn from the surface of a unit sphere. Then, we use $J(K_h; x_0)$ to denote an instantiation of the objective value $ J_h(K_h)$ given $x_0$, and define the two-point zeroth-order estimate of $\nabla  J_h^{r_h}$ as
\begin{align*}
	g(K_h) := \frac{mn}{2r_h}\big[J(K_h + r_hU; x_0) - J(K_h - r_hU; x_0)\big]U,
\end{align*}
where $mn$ is the dimension of the policy space. Subsequently, we invoke the $L_h$-smoothness property to derive
\begin{align*}
	&\hspace{1em}\EE^t[ J_h(K_{h, t+1}) -  J_h(K_{h, t})] \\
	&\leq \hspace{-0.1em}\EE^t \hspace{-0.1em}\Big[\hspace{-0.1em}\big\langle \nabla  J_h(K_{h, t}), K_{h, t+1}\hspace{-0.1em}-\hspace{-0.1em}K_{h, t}\big\rangle \hspace{-0.15em}+\hspace{-0.15em} \frac{L_h}{2}\hspace{-0.1em}\|K_{h, t+1} \hspace{-0.1em}-\hspace{-0.1em} K_{h, t}\|_F^2\hspace{-0.1em}\Big]\\
	&= \hspace{-0.1em}- \big\langle \eta_h\nabla J_h(K_{h, t}), \nabla J_h^{r_h}(K_{h, t})\hspace{-0.1em}\big\rangle \hspace{-0.1em}+\hspace{-0.1em} \frac{L_h\eta_h^2}{2}\EE^t\big[\|g(K_{h, t})\|_F^2\big]\\
	&= - \eta_h\|\nabla  J_h(K_{h, t})\|_F^2 + \eta_h L_h r_h \|\nabla J_h(K_{h, t})\|_F \\
	&\hspace{1em}+ \frac{L_h\eta_h^2}{2}\EE^t\big[\|g(K_{h, t})\|_F^2\big],
\end{align*}
where the inequalities are due to Lemma 14 of \cite{malik2020derivative}. Moreover, 
\begin{align*}
	\EE^t\big[\|g(K_{h, t})\|_F^2\big] &= \Var(g(K_{h, t})) + \|\nabla  J^{r_h}_h(K_{h, t})\|_F^2 \\
	&\leq \Var(g(K_{h, t})) + 2 \|\nabla  J_h(K_{h, t})\|_F^2 \\
	&\hspace{1em}+ 2 \|\nabla  J^{r_h}_h(K_{h, t})- \nabla  J_h(K_{h, t})\|^2_F \\
	&\leq G_2 + 2\|\nabla  J_h(K_{h, t})\|_F^2 + 2L_h^2 r_h^2.
\end{align*}
Again by the $L_h$-smoothness property, we have
\begin{align*}
	&\hspace{1em} J_h(K_{h, t} - \eta_h\nabla J_h(K_{h, t})) \\
	&\leq  J_h(K_{h, t}) - (\eta_h - \frac{\eta_h^2L_h}{2})\|\nabla J_h(K_{h, t})\|_F^2 \\
	&\hspace{-1em}\Longrightarrow (\eta_h - \frac{\eta_h^2L_h}{2})\|\nabla J_h(K_{h, t})\|_F^2 \\
	&\leq  J_h(K_{h, t}) -  J_h(K_{h, t} - \eta_h\nabla J_h(K_{h, t})) \\
	&\leq  J_h(K_{h, t}) -  J_h(\tilde{K}_h^*) = \Delta_t.
\end{align*}
Then, letting $\eta_h \in (0,  \frac{1}{2L_h}]$, we can derive
\begin{align*}
	&\hspace{1em}\EE^t[\Delta_{t+1} - \Delta_t] \\
	&\leq - \eta_h\|\nabla  J_h(K_{h, t})\|_F^2 + 2\eta_h L_h r_h \Delta_t^{1/2} + \frac{L_h\eta_h^2}{2}G_2 \\
	&\hspace{1em}+ L_h\eta_h^2\|\nabla  J_h(K_{h, t})\|_F^2 + \eta_h^2L_h^3 r_h^2 \\
	&\hspace{-0.2em}\leq \hspace{-0.1em}-\frac{\eta_h\alpha_h}{2}\Delta_t \hspace{-0.1em}+\hspace{-0.1em} \frac{\eta_h\alpha_h}{4}\hspace{-0.1em}\Delta_t \hspace{-0.1em}+\hspace{-0.1em} \frac{4\eta_h L_h^2 r_h^2}{\alpha_h} \hspace{-0.1em}+\hspace{-0.1em} \frac{L_h\eta_h^2}{2}G_2 + \eta_h^2L_h^3 r_h^2,
\end{align*}
where the second inequality is due to that the $\alpha_h$ strong-convexity implies the $\alpha_h$ gradient domination property, the choice of stepsize $\eta_h \leq \frac{1}{2L_h}$, and $2ab \leq a^2+b^2$ for any $a, b$. Recall the choices of algorithmic parameters as follows:
\begin{align*}
	\eta_h \leq \frac{1}{2L_h}, \quad r_h \leq \min\Big\{\frac{\alpha_h}{4L_h}\sqrt{\frac{\varsigma\delta}{10}}, \frac{1}{2L_h}\sqrt{\frac{\alpha_h\varsigma\delta}{5}}\Big\}.
\end{align*}
Then, using the bounds on algorithmic parameters and rearranging terms lead to 
\begin{align*}
	\EE^t[\Delta_{t+1}] \leq \Big(1-\frac{\eta_h\alpha_h}{4}\Big)\Delta_t + \frac{L_h\eta_h^2}{2}G_2 + \eta_h\alpha_h\frac{\varsigma\delta}{20},
\end{align*}
which completes the proof.

\end{document}